\documentclass[a4paper,10pt,reqno]{amsart}
\usepackage{graphicx}
\usepackage{geometry}
\usepackage{float}
\usepackage{verbatim}
\usepackage{setspace}
\usepackage{siunitx}
\usepackage{paralist}
\usepackage{enumitem}
\usepackage{xcolor}
\usepackage{amsfonts,amsmath,amssymb,bbm,amsthm,amsbsy}
\usepackage[utf8]{inputenc}
\usepackage{lipsum}
\usepackage[english]{babel}
\usepackage{booktabs}
\usepackage{array}
\usepackage{subfig}
\usepackage{caption}
\usepackage{mathtools}
\usepackage{stmaryrd}

\newenvironment{frcseries}{\fontfamily{frc}\selectfont}{}

\makeatletter
\newtheorem*{rep@theorem}{\rep@title}
\newcommand{\newreptheorem}[2]{
\newenvironment{rep#1}[1]{
 \def\rep@title{#2 \ref{##1}}
 \begin{rep@theorem}}
 {\end{rep@theorem}}}
\makeatother

\newtheorem{thm}{Theorem}[section]
\newtheorem{lemma}[thm]{Lemma}
\newtheorem{prop}[thm]{Proposition}
\newtheorem{corr}[thm]{Corollary}

\newtheorem*{thm*}{Theorem}
\newtheorem*{lemma*}{Lemma}
\newtheorem*{prop*}{Proposition}
\newtheorem*{corr*}{Corrolary}
\newtheorem*{claim*}{Claim}

\theoremstyle{remark}

\newtheorem*{rmk*}{Remark}
\newtheorem*{conj*}{Conjecture}
\newtheorem*{quest*}{Question}

\theoremstyle{definition}
\newtheorem{defn}[thm]{Definition}
\newtheorem{exmp}[thm]{Example}

\newtheorem*{defn*}{Definition}
\newtheorem*{exmp*}{Example}

\newreptheorem{theorem}{Theorem}
\newreptheorem{corollary}{Corollary}
\newreptheorem{proposition}{Proposition}

\usepackage{fancyhdr}
\usepackage{tikz}
\usepackage{tikz-cd}

\usepackage{IEEEtrantools}

\newenvironment{equ*}[1]{\begin{IEEEeqnarray*}{#1}}{\end{IEEEeqnarray*}}

\newcommand{\Z}{\mathbb{Z}}

\newcommand{\Qc}{\mathcal{Q}}

\newcommand{\Nc}{\mathcal{N}}

\newcommand{\K}{\mathcal{K}}

\newcommand{\X}{\mathcal{X}}
\newcommand{\Y}{\mathcal{Y}}

\DeclareFontFamily{U}{mathx}{}
\DeclareFontShape{U}{mathx}{m}{n}{<-> mathx10}{}
\DeclareSymbolFont{mathx}{U}{mathx}{m}{n}
\DeclareMathAccent{\widecheck}{0}{mathx}{"71}

\newcommand{\inj}{\hookrightarrow}
\newcommand{\sur}{\twoheadrightarrow}

\DeclareMathOperator{\Del}{\Delta}

\DeclareMathOperator{\botimes}{\otimes^\blacksquare}

\DeclareMathOperator{\Hom}{Hom}

\DeclareMathOperator{\UHom}{\underline{\textup{Hom}}}

\DeclareMathOperator{\Ind}{Ind}
\DeclareMathOperator{\Coind}{Coind}
\DeclareMathOperator{\Res}{Res}

\newcommand{\Set}{\mathbf{Set}}
\newcommand{\Grp}{\mathbf{Grp}}

\newcommand{\Ring}{\mathbf{Ring}}

\newcommand{\Top}{\mathbf{Top}}

\newcommand{\Pro}{\mathbf{Pro}}
\newcommand{\CHED}{\mathbf{CHED}}

\newcommand{\CSet}{\mathbf{CondSet}}
\newcommand{\CGrp}{\mathbf{CondGrp}}
\newcommand{\CAb}{\mathbf{CondAb}}

\newcommand{\CModRR}{\mathbf{CMod}(\mathcal{R})}

\newcommand{\SR}{\mathbf{Solid}(\mathcal{R})}

\newcommand{\R}{\mathcal{R}}
\newcommand{\M}{\mathcal{M}}
\newcommand{\G}{\mathcal{G}}
\newcommand{\Hc}{\mathcal{H}}

\title{Open Condensed Subgroups and Mackey's Formula}
\author{Jiacheng Tang}
\thanks{Email: jiacheng.tang@postgrad.manchester.ac.uk, University of Manchester}

\calclayout
\begin{document}
\maketitle

\begin{abstract}
We define what it means for a condensed group action to be open (following \cite{condopen}) and show that for open subgroups, many elementary results about abstract modules hold for condensed modules, such as the existence of Mackey's Formula for condensed groups. We also indicate how these results can be ``solidified" to obtain their solid versions.
\end{abstract}

\section{Introduction}
\label{sec1}

Condensed mathematics was recently developed as a categorical framework to combine algebra and topology, in a way that behaves better than the classical theory of topological groups and modules (see \cite{condensed}). One natural question to ask is how condensed modules and cohomology compare to other theories, such as abstract modules and cohomology. There are obvious ways we can define (co)induction and restriction (Definition \ref{inddef}), but it turns out that unlike the abstract case, restriction of condensed groups does not always preserve projectives (Lemma \ref{brinkeg}). Since this is such an important property of abstract modules, one is led to think about when restriction of condensed groups does preserve projectives.

In \cite{condopen}, Peter Scholze defines what it means for a condensed subgroup $\Hc$ of $\G$ to be \emph{open}, which happens if and only if the quotient condensed set $\G/\Hc$ comes from a discrete topological space (Proposition \ref{Gdecomp}). In this case, we can decompose $\G$ into a coproduct of cosets of $\Hc$, as in the abstract case, which is sufficient to prove many of the basic results about abstract modules or cohomology. For example, the paper \cite{zou} proves that restriction to an open subgroup does preserve projectives (\cite[Lemma 3.0.9]{zou}). As a further illustration, we will show that the condensed analogue of Mackey's Formula holds when we have an open subgroup (Theorem \ref{mackey}). The key takeaway from this paper should not be any of these results (since the proofs are all formal), but rather the fact that there is a good notion of openness, and that for open subgroups, many elementary results for abstract modules transfer directly to the condensed setting. This means that open subgroups are either interesting or uninteresting, depending on the reader's point of view.

We will assume knowledge of basic category theory (refer to \cite{maccat}), sheaf theory (refer to \cite{macsheaves}), condensed mathematics (refer to \cite{condensed}), and group cohomology (refer to \cite{brown}). Some of our examples will come from profinite groups (refer to \cite{profinite}), and the reader can find an introduction to condensed modules in the appendix of \cite{jc1}.

In Section \ref{sec2}, we will fix notations and define condensed group actions in the obvious way. In Section \ref{sec3}, we define (co)induction and restriction and discuss some of their basic properties. The main point, as mentioned above, is that restriction does not always preserves projectives (Lemma \ref{brinkeg}). In Section \ref{sec4}, we generalise the idea from \cite{condopen} slightly and define what it means for a condensed group action on a condensed set $\X$ to be open, which allows us to decompose $\X$ into its orbits (Propositions \ref{Gdecomp} and \ref{Gbide}). As a consequence, the condensed analogues of many results from abstract module theory and cohomology hold for open subgroups, as we shall see in Section \ref{sec5}, including Mackey's Formula (Theorem \ref{mackey}). Some of our corollaries were already proven in \cite{zou} when the condensed ring is $\underline{\Z}$, but its arguments clearly extend to all condensed rings. We point out that unlike the abstract case, the category $\CAb$ of condensed abelian groups does not have enough injectives (see \cite{condinj}). In Section \ref{sec6}, we indicate how some of the results from Section \ref{sec5} can be ``solidified" to obtain analogous statements for solid modules.

Remark: Whenever we write ``=" in this paper, we mean isomorphic, usually canonically isomorphic (or equivalent in the case of categories).

Convention: All rings are associative with a 1 but are not necessarily commutative.

\subsection*{Acknowledgements}
The author would like to thank his supervisor Peter Symonds for his constant guidance, the following people (in alphabetical order) for helpful discussions on the subject matter: Matthew Antrobus, Calum Hughes, and Peter Scholze, as well as the following person for reading drafts of this paper and giving useful feedback: Gregory Kendall.


\section{Notations}
\label{sec2}

Let $\CHED$ denote the category of compact Hausdorff extremally disconnected topological spaces. Recall the following definition from \cite{condensed} (ignoring set-theoretic issues):

\begin{defn}\label{condset}
A \emph{condensed set/group/ring/\ldots} is a sheaf of sets/groups/rings/\ldots on the site $\CHED$, with finite jointly surjective families of maps as covers. That is, a condensed set/group/ring/\ldots is a functor $$\mathcal{T}\colon\CHED^{\text{op}}\to\Set/\Grp/\Ring/\ldots$$ such that $\mathcal{T}(\varnothing)=*$ and for any $S_1, S_2\in\CHED$, the natural map $\mathcal{T}(S_1\sqcup S_2)\to\mathcal{T}(S_1)\times\mathcal{T}(S_2)$ is a bijection.

Given a condensed set/group/ring/\ldots $\mathcal{T}$, we call the global section $\mathcal{T}(*)$ its \emph{underlying set/\linebreak group/ring/\ldots}.
\end{defn}

Let $\mathbf{CondSet}$, $\mathbf{CondGrp}$, $\mathbf{CondAb}$ and $\mathbf{CondRing}$ denote respectively the (large) categories of condensed sets, condensed groups, condensed abelian groups, and condensed rings. Note that $\CSet$ and $\CGrp$ have all (small) limits and colimits. Limits and filtered colimits are computed pointwise, while general colimits have to be further sheafified. Given $T$ a ($T1$) topological space/group/ring/\ldots, we write $\underline{T}=C(-,T)$ for the associated condensed set/group/ring/\ldots.


\begin{defn}\label{condactiondef}
Let $\G$ be a condensed group and $\X$ a condensed set. We say that \emph{$\G$ acts on $\X$ (on the left)}, or that $\X$ is a \emph{(left) $\G$-set}, if there is a natural transformation $\G\times\X\to\X$ such that for each $S\in\CHED$, the function of sets $\G(S)\times\X(S)\to\X(S)$ makes $\X(S)$ a $\G(S)$-set. Equivalently, $\X$ is a (left) $\G$-object in the monoidal category $(\CSet,\times)$. Let us denote the category of $\G$-sets (and $\G$-maps) by $\CSet(\G)$.

The \emph{orbit set} or \emph{quotient set} $\G\backslash\X$ of the action is the condensed set defined pointwise by $(\G\backslash\X)(S)=\G(S)\backslash\X(S)$, which does define a sheaf.

We define right actions similarly.

Let $\G, \Hc$ be condensed groups and $\X$ a condensed set. We say that $\X$ is a \emph{$\G$-$\Hc$-biset} if it is a left $\G$-set and a right $\Hc$-set and the actions commute. Equivalently, $\X$ is a $\G$-$\Hc$-bimodule object in $(\CSet,\times)$. We define orbit sets here similarly.
\end{defn}

Remark: We can canonically identify $\G$-$\Hc$-bisets with left $(\G\times\Hc^{\text{op}})$-sets. Very often, we will not indicate whether we are talking about left or right actions: this should either be clear from context, or it does not matter which side we consider.

\begin{defn}
Let $\G$ be a condensed group. A \emph{(condensed) subgroup} of $\G$ is a condensed group $\Hc$ together with an injective map $\Hc\inj\G$ of condensed groups, i.e.\ a natural transformation $\Hc\to\G$ such that for each $S\in\CHED$, the map $\Hc(S)\to\G(S)$ is an injective group homomorphism. In this case, we write $\Hc\leq\G$. The \emph{quotient set} $\G/\Hc$ is the orbit set where $\Hc$ acts on $\G$ on the right by multiplication.

\end{defn}

\section{(Co)induction and Restriction}
\label{sec3}

Given a condensed ring $\R$, let $\CModRR$ denote the (large) category of condensed $\R$-modules. Given $\M,\Nc\in\CModRR$ (which are appropriately sided), we can define their tensor product $\M\otimes_{\R}\Nc\in\CAb$ pointwise by $(\M\otimes_\R\Nc)(S)=\M(S)\otimes_{\R(S)}\Nc(S)$. There is also an enriched Hom $\UHom_\R(\M,\Nc)\in\CAb$ defined by $\UHom_\R(\M,\Nc)(S)=\Hom_\R(\Z[\underline{S}]\otimes\M,\Nc)$, where $\otimes$ means $\otimes_{\Z}=\otimes_{\underline{\Z}}$. See \cite[pages 25-27]{jc1} for more details.

For a condensed ring $\R$ and a condensed group $\G$, let $\R[\G]=\R\otimes\Z[\G]$ denote the corresponding \emph{condensed group ring}. Note that $\R[\G]$-modules are precisely $\R$-modules which are also $\G$-modules where the actions commute.

Remark: It seems like the vast majority of existing literature on group rings assumes the ring to be commutative, for good reasons. The author does not believe that anything in this paper requires commutativity of the ring $\R$, but for other applications this assumption might be necessary.

\begin{defn}\label{inddef}
Let $\R$ be a (fixed) condensed ring and $\Hc\leq\G$ be condensed groups. If $\M$ is a (condensed) right $\R[\Hc]$-module, its \emph{induced module (from $\Hc$ to $\G$)} is $$\Ind^{\G}_{\Hc}\M=\M\otimes_{\R[\Hc]}\R[\G],$$ which is a right $\R[\G]$-module.

On the other hand, its \emph{coinduced module (from $\Hc$ to $\G$)} is $$\Coind^{\G}_{\Hc}\M=\UHom_{\R[\Hc]}(\R[\G],\M),$$ which is also a right $\R[\G]$-module.

If $\Nc$ is an $\R[\G]$-module, its \emph{restricted module (from $\G$ to $\Hc$)}, denoted by $\Res^{\G}_{\Hc}\Nc$, is simply $\Nc$ viewed as an $\R[\Hc]$-module via the inclusion $\R[\Hc]\inj\R[\G]$. We will often abuse notation and write the restricted module as just $\Nc$.
\end{defn}

Recall that given $\R,\R'$ condensed rings, $\K$ a right $\R$-module, $\M$ an $\R$-$\R'$-bimodule, and $\mathcal{N}$ a right $\R'$-module, we have the Hom-tensor adjunction (see \cite[Proposition A.21]{jc1}) $$\Hom_{\R}(\mathcal{K},\UHom_{\R'}(\M,\mathcal{N}))=\Hom_{\R'}(\mathcal{K}\otimes_\R\M,\mathcal{N}).$$ As a special case, we have the following.

\begin{prop}[\cite{zou} Lemma 3.0.8]\label{indleftad}
Let $\R$ be a condensed ring and $\Hc\leq\G$ condensed groups. Then induction is left adjoint to and coinduction is right adjoint to the restriction functor $\mathbf{CMod}(\R[\G])\to\mathbf{CMod}(\R[\Hc])$. In particular, induction preserves projectives.
\begin{proof}
We use the adjunction stated just before the proposition.

For induction, let $\Nc$ be a right $\R[\G]$-module. Note that $\UHom_{\R[\G]}(\R[\G],\Nc)$, where $\R[\G]$ is viewed naturally as an $\R[\Hc]$-$\R[\G]$-bimodule, is just $\Nc$ viewed as a right $\R[\Hc]$-module.

For coinduction, let $\K$ be a right $\R[\G]$-module. Note that $\K\otimes_{\R[\G]}\R[\G]$, where $\R[\G]$ is viewed naturally as an $\R[\G]$-$\R[\Hc]$-bimodule, is just $\K$ viewed as a right $\R[\Hc]$-module.
\end{proof}
\end{prop}

\begin{lemma}\label{indexact}
Let $\R$ be a fixed condensed ring and let $\Hc\leq\G$ be condensed groups. Then the induction functor $\Ind^{\G}_{\Hc}(-)\colon\mathbf{CMod}(\R[\Hc])\to\mathbf{CMod}(\R[\G])$ is exact. In particular, $\R[\G]$ is flat over $\R[\Hc]$ i.e.\ the functor $(-)\otimes_{\R[\Hc]}\R[\G]\colon\mathbf{CMod}(\R[\Hc])\to\CAb$ is exact.
\begin{proof}
It suffices to prove exactness on the level of presheaves, where we can do so pointwise, but this is obvious because in the abstract case, $R[G]$ is free (so flat) over $R[H]$.
\end{proof}
\end{lemma}

As stated in the proof above, one useful phenomenon that occurs in the abstract (or profinite) case is that given groups $H\leq G$, the group ring $R[G]$ is free as an $R[H]$-module, which then implies that restriction from $G$ to $H$ preserves projectives. The analogue of this in the condensed setting is false:

\begin{lemma}[\cite{brink} Lemma A.2.1]\label{brinkeg}
Suppose $X$ is a (T1) compact topological space such that $\Z[\underline{X}]$ is a projective condensed abelian group. Then $X$ is totally disconnected.
\end{lemma}

Take any non-totally disconnected compact topological group $G$, such as the circle group $S^1$. Then $\Z[\underline{G}]=\underline{\Z}[\underline{G}]$ is not projective in $\CAb$ i.e.\ $\Z[\underline{G}]$ is not projective when restricted to the trivial subgroup. In particular, unlike induction (Lemma \ref{indexact}), the coinduction functor $\Coind^{\G}_{\Hc}(-)\colon\mathbf{CMod}(\R[\Hc])\to\mathbf{CMod}(\R[\G])$ isn't always exact.

One issue here is that, unlike in the abstract case, we cannot always express the condensed group $\G$ as a disjoint union of copies of the subgroup $\Hc$. However, we should be able to do so if $\Hc$ is ``open" in $\G$, in which case the cosets of $\Hc$ in $\G$ should form an ``open cover" of $\G$. We take inspiration from \cite{condopen} for the next definition.

\section{Open Actions}
\label{sec4}

\begin{defn}\label{openactdef}
Let $\G$ be a condensed group and $\X$ a left $\G$-set. Given $x\in\G\backslash\X(*)$, the \emph{orbit} $\G x\in\CSet$ is the pullback
\[\begin{tikzcd}
\G x \arrow[r] \arrow[d, hook] & \underline{*} \arrow[d, "x", hook]       \\
\X \arrow[r]                   & \G\backslash\X              
\end{tikzcd}\]

We say that the $\G$-action on $\X$ is \emph{open} if for any $S\in\CHED$, any map $\underline{S}\to\X$ and any $x\in\G\backslash\X(*)$, the condensed set pullback $\underline{S}\times_{\X}\G x$ is representable by an open subset of $S$. We define right orbits and open right actions similarly.

Given $\Hc\leq\G$ condensed group, we say that $\Hc$ is \emph{open} in $\G$ if for any $S\in\CHED$ and any map $\underline{S}\to\G$, the condensed set pullback $\underline{S}\times_{\G}\Hc$ is representable by an open subset of $S$.
\end{defn}

\begin{exmp}\label{openact}
\begin{enumerate}[label=(\roman*)]
\item Let $\Hc\leq\G$ be condensed groups. If the action of $\Hc$ on $\G$ (by right multiplication) is open, then clearly $\Hc$ is an open subgroup of $\G$ (pick $x=1\in\G/\Hc(*)$ in the second definition above). Conversely, if $\Hc$ is open in $\G$, then the right action of $\Hc$ on $\G$ is open. To see this, note that given $g\in\G/\Hc(*)$ and a map $\underline{S}\to\G$, we can view the pullback $\underline{S}\times_{\G}g\Hc$ as the pullback of the composite $\underline{S}\to\G\to\G$ along $\Hc\inj\G$, where the map $\G\to\G$ is ``multiplication by $g^{-1}$" (or more accurately, multiplication by $(g')^{-1}$, where $g'\in\G(*)$ is a lift of $g$).
\item\label{openact2} As pointed out in \cite{condopen}, if $G$ is a (T1) topological group and $H$ is an open (so clopen) subgroup, then $\underline{H}$ is open in $\underline{G}$ in the above sense. Indeed, given a map $\underline{S}\to\underline{G}$ i.e.\ a continuous map $f\colon S\to G$, let $P=f^{-1}(H)$, which is a clopen subset of $S$, so in particular also in $\CHED$. Then $\underline{P}$ is the required representable pullback since $\underline{(-)}$ preserves limits.
\item\label{topegopen} More generally, let $G$ be a topological group acting continuously on a topological space $X$, such that the quotient space $G\backslash X$ is discrete, and such that for each $x\in X$, the map $G\to Gx, g\mapsto gx$ has a continuous section (for example, these conditions are satisfied in (ii)). Then the induced action of $\underline{G}$ on $\underline{X}$ is open.
\item It is clear from the definition that a $\G$-action on $\X$ is open if and only if the trivial $\G$-action on $\G\backslash\X$ is open.
\end{enumerate}
\end{exmp}

Remark: In Definition \ref{openactdef}, we could have replaced $\CHED$ everywhere with $\Pro$, the category of profinite spaces, if we wanted. Indeed, suppose for any $S\in\CHED$, the pullback $\underline{S}\times_\X \G x$ is representable by an open subset of $S$. Given a profinite space $T$, we can find some $S\in\CHED$ with a surjection $S\sur T$, which is automatically a quotient map. Let $\Qc=\underline{T}\times_\X \G x$ and $\underline{P}=\underline{S}\times_{\underline{T}}\Qc=\underline{S}\times_\X \G x$, where $P\subseteq S$ is open. By giving $\Qc(*)\subseteq T$ the subspace topology from $T$, we see that $P=S\times_T\Qc(*)$ is a pullback in $\Top$ (and $\Qc(*)\subseteq T$ is clopen because $S\sur T$ is a quotient map), so that $\underline{P}=\underline{S}\times_{\underline{T}}\underline{\Qc(*)}$. Finally, the map $\underline{S}\to\underline{T}$ is an epimorphism, so $\Qc=\underline{\Qc(*)}$ is representable by an open subset of $T$.

Given a set/group/ring\ldots $A$, let $\Del(A)$ denote the constant presheaf on $\CHED$ with value $A$, and let $\Del_0(A)$ denote the presheaf defined by $\Del_0(A)(S)=A$ if $S\neq\varnothing$, and $\Del_0(A)(\varnothing)=*$. By mimicking the proof of \cite[Lemma A.4]{jc1}, we see that the sheafifications of $\Del(A)$ and $\Del_0(A)$ are both $\underline{A}=C(-,A)$, where $A$ is given the discrete topology.

\begin{lemma}\label{disclem}
Let $X$ be a (discrete) set and $\Y$ a condensed subset of $\underline{X}$ i.e.\ there is an injection $\Y\inj\underline{X}$. Then $\Y$ is discrete i.e.\ $\Y=\underline{\Y(*)}$, where $\Y(*)$ is given the discrete topology.
\begin{proof}
Let $Y=\Y(*)$. There is clearly an injection of presheaves $\Del_0(Y)\inj\Y$ and hence of sheaves $\underline{Y}\inj\Y$, so it suffices to show that the final map is pointwise surjective (this is in generally stronger than being a sheaf epimorphism). Given $S\in\CHED$ and an element $f\in\Y(S)$ i.e.\ a map $f\colon\underline{S}\to\Y$, its image $f'$ in $\underline{X}(S)=C(S,X)$ has finite image as a function (since $S$ is compact and $X$ is discrete), so it suffices to prove that $f$ has a preimage in $\underline{Y}(S)$ whenever $f'$ is a constant function, say with value $x\in X$. But then $x\in Y$, so we're done.
\end{proof}
\end{lemma}

\begin{lemma}\label{discalways}
Let $X$ be a (discrete) set and $\G$ a condensed group acting on $\underline{X}$. Then the action is open (cf.\ Example \ref{openact}\ref{topegopen}).
\begin{proof}
Suppose we are given $S\in\CHED$ and a map $\underline{S}\to\underline{X}$ i.e.\ a continuous map $S\to X$. For each $x\in\G(*)\backslash X$, the orbit $\G x\subseteq\underline{X}$ is discrete by Lemma \ref{disclem}, so we can simply take the pullback $S\times_X \G(*)x$ in $\Top$.
\end{proof}
\end{lemma}

We can easily obtain the following result by generalising the proof in \cite{condopen}.

\begin{prop}\label{Gdecomp}
Let $\G$ be a condensed group and $\X$ a $\G$-set. Then the action is open if and only if the quotient set $\G\backslash\X$ is discrete. In this case, we have an isomorphism $\X=\coprod_{x\in\G\backslash\X(*)}\G x$ of $\G$-sets, which decomposes $\X$ into its orbits.
\begin{proof}
Let $G=\G(*)$ and $X=\X(*)$ be the underlying group/set of the objects in question. If the quotient $\G\backslash\X$ is discrete i.e.\ $\G\backslash\X=\underline{G\backslash X}$, where $G\backslash X$ is given the discrete topology, then by Lemma \ref{discalways}, the trivial $\G$-action on $\G\backslash\X$ is open, so the original $\G$-action on $\X$ is open.

Conversely, assume that the $\G$-action on $\X$ is open. Suppose we are given $S\in\CHED$ and a map $\underline{S}\to\X$. For each $x\in G\backslash X$, let $S_x\subseteq S$ be the open subset representing the pullback $\underline{S}\times_{\X}\G x$. By looking at the underlying sets, we see that $S=\coprod_{G\backslash X}S_x$, where only finitely many of the $S_x$ are non-empty by compactness. Thus, the family $\{S_x\inj S\}_{G\backslash X}$ covers $S$ in the Grothendieck topology on $\CHED$. Similarly, given a map $\underline{S}\to\G\backslash\X$ and $x\in G\backslash X$, the pullback of $\underline{S}\to\G\backslash\X$ along $x\colon\underline{*}\to\G\backslash\X$ is also representable by an open subset $S_x'\subseteq S$, and these form a covering family $\{S_x'\inj S\}_{G\backslash X}$.

Recall that $\Del_0(A)$ is the presheaf on $\CHED$ that takes value $*$ at $\varnothing$ and value $A$ elsewhere. There is an injection of presheaves (of sets) $\phi\colon\Del_0(G\backslash X)\inj\G\backslash\X$ and hence an injection of sheaves $\underline{G\backslash X}\inj\G\backslash\X$. The covering family $\{S_x'\inj S\}_{G\backslash X}$ defined above shows that $\phi$ is locally surjective, so that $\underline{G\backslash X}=\G\backslash\X$ (see \cite[Corollary III.7.6]{macsheaves}) i.e.\ $\G\backslash\X$ is discrete. Similarly, we have an injection of presheaves and hence of sheaves $\coprod_{G\backslash X}\G x\inj\X$. The covering family $\{S_x\inj S\}_{G\backslash X}$ from above then shows this is an isomorphism. Note that this is not just an isomorphism of condensed sets, but also of $\G$-sets with the obvious actions. Indeed, since sheafification commutes with finite products and arbitrary coproducts, it suffices to check this on the level of presheaves, where it's clear.
\end{proof}
\end{prop}

\begin{exmp}\label{deceg}
\begin{enumerate}[label=(\roman*)]
\item\label{deceg1} In the special case of the above proposition when $\Hc$ is an open subgroup of a condensed group $\G$, we obtain the decomposition $\G=\coprod_{\G/\Hc(*)}g\Hc$ of \cite{condopen}. Note that a condensed subgroup $\Hc$ is open in $\G$ if and only if the quotient set $\G/\Hc$ is discrete.
\item If $G$ is an abstract group acting on a set $X$, we always have an orbit decomposition $X=\coprod_{G\backslash X}Gx$. The functor $\underline{(-)}\colon\Set\to\CSet$ preserves colimits (being left adjoint to the global sections functor), so we get the condensed set decomposition $\underline{X}=\coprod_{G\backslash X}\underline{G}x$. Indeed, by Lemma \ref{discalways}, the action of $\underline{G}$ on $\underline{X}$ is open, so the above proposition generalises orbit decompositions of discrete group actions.
\item Let $G$ be a profinite group acting continuously on a profinite space $X$. Then we cannot in general decompose $X$ into its (profinite) orbits, but if the conditions of Example \ref{openact}\ref{topegopen} are satisfied (in particular for the multiplication action by an open subgroup), then we do have a decomposition $X=\coprod_{G\backslash X}Gx$ into a finite disjoint union, and applying $\underline{(-)}$ shows that the above proposition generalises orbit decompositions of some profinite group actions.
\end{enumerate}
\end{exmp}

We can actually generalise the claims of Proposition \ref{Gdecomp} into ones for bisets, which we will need later.

\begin{prop}\label{Gbide}
Let $\G, \Hc$ be condensed groups and $\X$ a $\G$-$\Hc$-biset such that the left $\G$-action on $\X$ is open (or the right $\Hc$-action is open). Then the quotient set $\G\backslash\X/\Hc$ is discrete. Moreover, we have an isomorphism $\X=\coprod_{x\in\G\backslash\X/\Hc(*)}\G x\Hc$ of $\G$-$\Hc$-bisets, which decomposes $\X$ into its double orbits.
\begin{proof}
Let $G=\G(*)$, $H=\Hc(*)$ and $X=\X(*)$. We can view $\G\backslash\X/\Hc$ as the quotient set of the induced right $\Hc$-action on $\G\backslash\X$, so $\G\backslash\X/\Hc$ is discrete by Lemma \ref{discalways} and Proposition \ref{Gdecomp}.

Although the following proof might appear confusing, all we are doing is decomposing $\X$ into its left $\G$-orbits, and then grouping them into $\G$-$\Hc$-double orbits. Choose set sections $f$ of $X\sur G\backslash X$ and $g$ of $G\backslash X\sur G\backslash X/H$. By Proposition \ref{Gdecomp} again, we have the orbit decomposition $\X=\coprod_{y\in G\backslash X}\G f(y)$. As underlying sets, there is certainly a decomposition $G\backslash X=\coprod_{x\in G\backslash X/H}g(x)H$, so we have $\X=\coprod_{x\in G\backslash X/H}\coprod_{y\in g(x)H}\G f(y)$.

On the other hand, for each $x\in G\backslash X/H$, there is an induced action of $\G$ on the orbit $\G [fg(x)]\Hc\subseteq\X$ which is still an open action. We can canonically identify the underlying quotient set $G\backslash(G[fg(x)]H)$ with $g(x)H$ as subsets of $G\backslash X$, so by Proposition \ref{Gdecomp}, we have $$\X=\coprod_{x\in G\backslash X/H}\coprod_{y\in g(x)H}\G f(y)=\coprod_{x\in G\backslash X/H}\G[fg(x)]\Hc.$$

Although we have chosen some set sections to improve readability, the final decomposition is canonical and does not depend on any lifts.
\end{proof}
\end{prop}

Remark: An alternative proof of the above proposition is to directly mimic Proposition \ref{Gdecomp}. We shall leave it to the reader to define what it means for a $\G$-$\Hc$-action to be \emph{open} (mimic Definition \ref{openactdef}), but as suggested by the above proposition, this is implied by either the $\G$-action or the $\Hc$-action being open, so is not that interesting.

\section{Corollaries and Mackey's Formula}
\label{sec5}

We shall now see that for open subgroups, (co)induction and restriction of condensed groups behave very much like in the abstract case (apart from anything that involves injective modules; see \cite{condinj}). The author later discovered that some of these have already been proven in \cite{zou} when the condensed ring is $\underline{\Z}$, but its arguments clearly work for any condensed ring.

\begin{corr}[\cite{zou} Lemma 3.0.9]\label{preproj}
Let $\R$ be a fixed condensed ring and let $\Hc\leq\G$ be condensed groups with $\Hc$ open. Then restriction from $\G$ to $\Hc$ preserves projectives.
\begin{proof}
It suffices to show that restriction preserves the projectivity of the projective generators $\Z[\underline{S}]\otimes\R[\G]$, $S\in\CHED$ (see \cite[Theorem A.15]{jc1}). Note that the functor $\Z[-]\colon(\CSet,\times)\to(\CAb,\otimes)$ is monoidal (see \cite[Corollary A.11]{jc1}), so induces a functor $\CSet(\Hc)\to\mathbf{CMod}(\Z[\Hc])$ which is left adjoint to the forgetful functor $\mathbf{CMod}(\Z[\Hc])\to\CSet(\Hc)$. Thus, by applying $\Z[-]$ to the decomposition of $\G$ in Example \ref{deceg}\ref{deceg1}, we see that $\Z[\G]=\bigoplus_{\G/\Hc(*)}\Z[\Hc]$ as $\Z[\Hc]$-modules. Now note that the forgetful functor $\mathbf{CMod}(\R[\Hc])\to\mathbf{CMod}(\Z[\Hc])$ has a left adjoint $(-)\otimes\R\colon\mathbf{CMod}(\Z[\Hc])\to\mathbf{CMod}(\R[\Hc])$, so $\Z[\underline{S}]\otimes\R[\G]=\bigoplus_{\G/\Hc(*)}(\Z[\underline{S}]\otimes\R[\Hc])$ as right $\R[\Hc]$-modules for each $S\in\CHED$, completing the proof.
\end{proof}
\end{corr}

Remark: Note that \cite[Lemma 3.0.9]{zou} proves this by establishing a non-canonical $\Hc$-isomorphism $\G=(\G/\Hc)\times\Hc$, whereas the decomposition $\G=\coprod_{\G/\Hc(*)}g\Hc$ from Example \ref{deceg}\ref{deceg1} is canonical and does not depend on a section $\G/\Hc(*)\to\G(*)$.

One important notion in abstract group theory is the property of having finite index: a subgroup $H$ has finite index in $G$ if the set $G/H$ is finite. But what does ``finite" mean in the condensed setting? Recall that to avoid set-theoretic difficulties, we only want to consider $T1$ topological spaces (see \cite[Proposition 2.15]{condensed}). Now a finite $T1$ topological space is discrete, so a condensed set should be called \emph{finite} precisely when it is discrete and the underlying set is finite. This seems like a natural way to define ``finite index" in the condensed world.

\begin{defn}
Let $\Hc\leq\G$ be condensed groups. We say that $\Hc$ \emph{has finite index in $\G$} if $\Hc$ is open in $\G$ and the quotient $\G/\Hc(*)$ is finite i.e.\ if the $\G/\Hc$ is finite as a condensed set.
\end{defn}

\begin{corr}[\cite{zou} Lemma 3.0.11]\label{indiscoind}
Let $\R$ be a fixed condensed ring and let $\Hc\leq\G$ be condensed groups with $\Hc$ of finite index. Then induction and coinduction (from $\Hc$ to $\G$) are canonically isomorphic.
\begin{proof}
By Proposition \ref{Gdecomp} and the argument in Corollary \ref{preproj}, we have an $\R[\Hc]$-isomorphism $\R[\G]=\bigoplus_{\G/\Hc(*)}\R[\Hc]$, where the direct sum is finite. We then have, for $\M\in\mathbf{CMod}(\R[\Hc])$, that $$\Ind^\G_{\Hc}\M=\bigoplus_{\G/\Hc(*)}\M=\prod_{\G/\Hc(*)}\M=\Coind^\G_{\Hc}\M.$$ This might seem to only be an isomorphism of $\R[\Hc]$-modules (or even of condensed abelian groups), but it is in fact canonically an $\R[\G]$-module isomorphism, as we will explain below.
\end{proof}
\end{corr}

Let us now explain why the isomorphism in the proof above is an $\R[\G]$-module isomorphism. Recall the decomposition $\G=\coprod_{\G/\Hc(*)}\Hc g$ of Example  \ref{deceg}\ref{deceg1}, where $\Hc$ is an open subgroup of $\G$ (now acting on $\G$ on the left). We only claimed that this is an isomorphism of left $\Hc$-sets, but in fact we can use this isomorphism (of condensed sets) to induce a $\G$-$\G$-action on $\coprod_{\G/\Hc(*)}\Hc g$, exactly as in the abstract case. Note that this new left $\G$-action is compatible with the original left $\Hc$-action. In this way, the decomposition (tautologically) becomes a $\G$-$\G$-biset isomorphism. Applying $\R[-]$ then gives an $\R[\G]$-$\R[\G]$-bimodule isomorphism $\R[\G]=\bigoplus_{\G/\Hc(*)}\R[\Hc g]$, which is what we used in the proof above.

We can play a similar game with double cosets. Let $\Hc,\K\leq\G$ be condensed groups with $\Hc$ open (or $\K$ open). Then by Proposition \ref{Gbide}, we have an isomorphism $\G=\coprod_{\Hc\backslash\G/\K(*)}\Hc g\K$ of $\Hc$-$\K$-bisets. Applying $\R[-]$ then gives an $\R[\Hc]$-$\R[\K]$-bimodule isomorphism $$\R[\G]=\bigoplus_{\Hc\backslash\G/\K(*)}\R[\Hc g\K].$$ If $\M\in\mathbf{CMod}(\R[\Hc])$, then applying $\M\otimes_{\R[\Hc]}(-)$ to the above gives an $\R[\K]$-isomorphism $$\Res^\G_\K\Ind^\G_\Hc\M=\M\otimes_{\R[\Hc]}\R[\G]=\bigoplus_{\Hc\backslash\G/\K(*)}\M\otimes_{\R[\Hc]}\R[\Hc g\K].$$ Hopefully, this is starting to look familiar.

Given two subgroups $\Hc,\K\leq\G$, we can define their \emph{intersection} $\Hc\cap\K$ pointwise, which is a condensed subgroup of $\G$. Alternatively, it is the pullback $\Hc\times_\G\K$. Given $\Hc\leq\G$ and $g\in\G(*)$, we can naturally form the \emph{conjugate} $g^{-1}\Hc g$, which is a condensed subgroup of $\G$. It is easy to check that if $\Hc$ and $\K$ are open subgroups of $\G$, then so are $\Hc\cap\K$ and $g^{-1}\Hc g$. The condensed set $\Hc g$ is a right $g^{-1}\Hc g$-set, so we can restrict it to become a $\K\cap g^{-1}\Hc g$-set.

\begin{lemma}
Let $\Hc,\K\leq\G$ be condensed groups and $g\in\G(*)$. Then there is an $\R[\Hc]$-$\R[\K]$-bimodule isomorphism $\R[\Hc g]\otimes_{\R[\K\cap g^{-1}\Hc g]}\R[\K]=\R[\Hc g\K]$.
\begin{proof}
It suffices to show this on the level of presheaves, but we know this is true in the abstract case. Explicitly, the (presheaf) isomorphism is given left-to-right by $hg\otimes k\mapsto hgk$.
\end{proof}
\end{lemma}

Combining what we have above, we finally obtain:

\begin{thm}[Mackey's Formula]\label{mackey}
Let $\Hc,\K\leq\G$ be condensed groups with $\Hc$ open (or $\K$ open) and let $\M$ be a right $\R[\Hc]$-module, where $\R$ is a fixed condensed ring. Then there is an $\R[\K]$-module isomorphism $$\Res^\G_\K\Ind^\G_\Hc\M=\bigoplus_{g\in\Hc\backslash\G/\K(*)}\Ind^\K_{\K\cap g^{-1}\Hc g}\Res^{g^{-1}\Hc g}_{\K\cap g^{-1}\Hc g}(\M\otimes_{\R[\Hc]}\R[\Hc g]).$$
\end{thm}

\section{The Solid Theory}
\label{sec6}

We shall now consider solid modules (refer to \cite[Lecture V]{condensed} or \cite[Section 3.2]{jc1}). For a condensed ring $\R$, let $\SR$ denote the category of solid $\R$-modules i.e.\ condensed $\R$-modules which are solid as condensed abelian groups. Given a condensed ring $\R$ and a condensed group $\G$, their \emph{solid group ring} is the solid ring $\R[\G]^\blacksquare=\R^\blacksquare\botimes\Z[\G]^\blacksquare$. Note that there is a canonical equivalence between the category of solid $\R[\G]^\blacksquare$-modules and the category of solid $\R[\G]$-modules.

We can define the solid versions of most of the concepts above in the obvious way. For example, if $\Hc\leq\G$ is a condensed subgroup and $\M\in\mathbf{Solid}(\R[\Hc]^\blacksquare)$, its \emph{induced (solid) module} (from $\Hc$ to $\G$) is $$\Ind^{\G\blacksquare}_\Hc\M=(\Ind^\G_\Hc\M)^\blacksquare=\M\otimes^\blacksquare_{\R[\Hc]^\blacksquare}\R[\G]^\blacksquare,$$ while its \emph{coinduced module} is simply $$\Coind^{\G\blacksquare}_\Hc\M=\Coind^\G_\Hc\M=\UHom_{\R[\Hc]^\blacksquare}(\R[\G]^\blacksquare,\M),$$ since this is already solid (see \cite[Lemma 3.13]{jc1}). One can easily state and prove the solid analogues of Proposition \ref{indleftad}, Corollary \ref{preproj}, Corollary \ref{indiscoind} and Theorem \ref{mackey} by simply solidifying the corresponding result/argument from above (noting that the solidification functor preserves colimits and is monoidal etc., see \cite[Theorems 5.8 and 6.2]{condensed} or \cite[page 13]{jc1}). Let us do this for Theorem \ref{mackey} (Mackey's Formula) as an illustration.

Let $\Hc,\K\leq\G$ be condensed groups with $\Hc$ open and let $\M$ be a right solid $\R[\Hc]^\blacksquare$-module, where $\R$ is a fixed condensed ring. By treating $\M$ as a (solid) $\R[\Hc]$-module, we obtain from Theorem \ref{mackey} an $\R[\K]$-module isomorphism $$\Res^\G_\K\Ind^\G_\Hc\M=\bigoplus_{g\in\Hc\backslash\G/\K(*)}\Ind^\K_{\K\cap g^{-1}\Hc g}\Res^{g^{-1}\Hc g}_{\K\cap g^{-1}\Hc g}(\M\otimes_{\R[\Hc]}\R[\Hc g]).$$ Now we can apply the solidification functor $(-)^\blacksquare$ to the above formula to get a (solid) $\R[\K]^\blacksquare$-module isomorphism $$\Res^{\G\blacksquare}_\K\Ind^{\G\blacksquare}_\Hc\M=\bigoplus_{g\in\Hc\backslash\G/\K(*)}\Ind^{\K\blacksquare}_{\K\cap g^{-1}\Hc g}\Res^{g^{-1}\Hc g\blacksquare}_{\K\cap g^{-1}\Hc g}(\M\otimes^\blacksquare_{\R[\Hc]^\blacksquare}\R[\Hc g]^\blacksquare),$$ which is the solid version of Mackey's Formula.

Note that this is a direct generalisation of Mackey's Formula for profinite groups (\cite[Proposition 6.11.2]{profinite}). Indeed, if $R$ is a profinite ring and $H,K\leq G$ are profinite groups with $H$ open, then the condensation functor $\underline{(-)}$ from profinite modules to solid modules preserves all operations appearing in Mackey's Formula. To be precise, note that $\underline{H}$ is open in $\underline{G}$ (Example \ref{openact}\ref{openact2}), and that $\underline{(-)}$ preserves groups rings (\cite[Lemma B.4(i)]{guido}), tensor products (\cite[Proposition 3.20(iii)]{jc1}), and limits (\cite[Lemma 3.8(i)]{jc1}), so in particular finite coproducts and pullbacks.

We point out that the solid analogue of Lemma \ref{indexact} is not necessarily true. In fact, the author believes that it is probably false. An unpublished example of Efimov shows that the solid abelian group $\prod_I\underline{\Z}$ is not flat with respect to $\botimes$ when $|I|=2^{2^{\aleph_0}}$, so if there exists a profinite group $G$ such that $C(G,\Z)=\bigoplus_I\Z$ for this index set $I$, then $\Z[\underline{G}]^\blacksquare=\prod_I{\underline{\Z}}$ is not flat over $\underline{\Z}$.

Let us take a look at the solid analogue of Corollary \ref{preproj}, which says that for $\Hc\leq\G$ an open subgroup, (solid) restriction preserves projectives. Suppose $\Hc=\underline{H}$ and $\G=\underline{G}$ come from profinite groups $H\leq G$, and the condensed ring is also profinite. Then we know that restriction \emph{always} preserves projective profinite modules, regardless of whether $H$ is open (see \cite[Corollary 5.7.2(b)]{profinite}), so it would seem like the solid analogue of Corollary \ref{preproj} is a lot weaker than what it should be. Indeed, for profinite groups we can strengthen the result.

\begin{lemma}
Let $\R$ be a condensed ring (not necessarily profinite) and $H\leq G$ profinite groups. Then the restriction functor $\mathbf{Solid}(\R[\underline{G}]^\blacksquare)\to\mathbf{Solid}(\R[\underline{H}]^\blacksquare)$ preserves projectives.
\begin{proof}
As profinite left $H$-sets, we have (non-canonically) $G=H\times G/H$, so $\underline{G}=\underline{H}\times\underline{G/H}$. Applying $\R[-]^\blacksquare$ shows that $\R[\underline{G}]^\blacksquare=\R[\underline{H}]^\blacksquare\otimes^\blacksquare\Z[\underline{G/H}]^\blacksquare$ as $\R[\underline{H}]^\blacksquare$-modules, so for each $S\in\CHED$, we have $\R[\underline{G}]^\blacksquare\botimes\Z[\underline{S}]^\blacksquare=\R[\underline{H}]^\blacksquare\otimes^\blacksquare\Z[\underline{(G/H)\times S}]^\blacksquare$, as required. Note that we used the fact that $\Z[\underline{T}]^\blacksquare$ is a projective solid abelian group for every profinite space $T$, rather than just for the extremally disconnected ones.
\end{proof}
\end{lemma}

Remark: The same argument shows that more generally, for $\R$ a condensed ring and $\Hc\leq\G$ condensed groups, restriction of solid modules from $\G$ to $\Hc$ preserves projectives if the quotient $\G/\Hc$ is representable i.e.\ if $\G/\Hc=\underline{X}$ for some profinite space $X$. (Here, we are viewing condensed objects as sheaves on $\Pro$ rather than $\CHED$.) Note that this observation is not entirely trivial: we might not be able to write $\G=\Hc\times\G/\Hc$, since there might not be a (natural) section of $\G\sur\G/\Hc$. However, if $\G/\Hc=\underline{X}$ is profinite, then to give a map $\underline{X}\to\G$ is the same as to give an element of $\G(X)$, which we can choose to be any preimage under $\G(X)\sur\underline{X}(X)$ of the identity map on $X$. This will be a section of $\G\sur\underline{X}$ and gives us the required (non-canonical) decomposition $\G=\Hc\times\G/\Hc$.

\bibliographystyle{unsrt}

\begin{thebibliography}{10}

\bibitem{condopen}
Peter Scholze.
\newblock ``{S}tructure of a profinite group as a condensed set with an action
  of an open subgroup", 2022.
\newblock URI:
  https://mathoverflow.net/questions/427868/structure-of-a-profinite-group-as-a-condensed-set-with-an-action-of-an-open-subg.

\bibitem{condensed}
Peter Scholze.
\newblock Condensed mathematics.
\newblock {\em Lecture notes based on joint work with D. Clausen. Available at
  this link https://www.math.uni-bonn.de/people/scholze/Condensed.pdf}, 2019.

\bibitem{zou}
Konrad Zou.
\newblock The categorical form of fargues' conjecture for tori.
\newblock {\em arXiv preprint arXiv:2202.13238}, 2022.

\bibitem{maccat}
Saunders Mac~Lane.
\newblock {\em Categories for the working mathematician}, volume~5.
\newblock Springer Science \& Business Media, 2013.

\bibitem{macsheaves}
Saunders MacLane and Ieke Moerdijk.
\newblock {\em Sheaves in geometry and logic: A first introduction to topos
  theory}.
\newblock Springer Science \& Business Media, 2012.

\bibitem{brown}
Kenneth~S. Brown.
\newblock {\em Cohomology of groups}, volume~87 of {\em Graduate Texts in
  Mathematics}.
\newblock Springer-Verlag, New York, 1994.
\newblock Corrected reprint of the 1982 original.

\bibitem{profinite}
Luis Ribes and Pavel Zalesskii.
\newblock {\em Profinite groups}.
\newblock Springer, 2000.

\bibitem{jc1}
Jiacheng Tang.
\newblock Profinite and solid cohomology.
\newblock {\em arXiv preprint arXiv:2410.08933}, 2024.

\bibitem{condinj}
Peter Scholze.
\newblock ``{A}re there (enough) injectives in condensed abelian groups?",
  2020.
\newblock URI:
  https://mathoverflow.net/questions/352448/are-there-enough-injectives-in-condensed-abelian-groups.

\bibitem{brink}
Emma Brink.
\newblock Condensed group cohomology.
\newblock Master's thesis, Ludwig-Maximilians-Universit\"at M\"unchen, 2023.

\bibitem{guido}
Guido Bosco.
\newblock On the $ p $-adic pro-\'etale cohomology of drinfeld symmetric
  spaces.
\newblock {\em arXiv preprint arXiv:2110.10683}, 2021.

\end{thebibliography}

\end{document}